\documentclass[11pt,titlepage,a4paper]{article}
\usepackage{bozhomac}	
\usepackage{bozhlogo}	
\usepackage{cite}	

\title{\bfseries	\vspace*{-1.678902345in}
{\huge Affine and fundamental vector fields}
}

\vspace{1.7ex}

\author{
Bozhidar Z.\ Iliev
\thanks{Laboratory of Mathematical Modeling in Physics,
Institute for Nuclear Research and \mbox{Nuclear} Energy,
Bulgarian Academy of Sciences,
Boul.\ Tzarigradsko chauss\'ee~72, 1784 Sofia, Bulgaria}
\thanks{E-mail address: bozho@inrne.bas.bg}
\thanks{URL: http://theo.inrne.bas.bg/$\sim$bozho/}
}

\date{
 \vspace{2.27ex}\ShortTitle{Affine and fundamental vector fields}	\\[0.27ex]
 \vspace{3.27ex}
\small
	\begin{tabular}{r@{$\colon\to~$}l}
\vspace{0.09ex} Last update	& February 1, 2006		\\[0.09ex]
 \vspace{0.27ex} Produced	& \fbox{\today}	\\[0.27ex]
	\end{tabular} \\[1.27ex]
 \small
	\begin{tabular}{r@{$\colon~$}l}
 \normalsize\sffamily\bfseries
  \vspace{0.27ex} http://www.arXiv.org e-Print archive No. &
 \normalsize\sffamily\bfseries
 math.DG/0602006 \\[1.27ex]
	\end{tabular} \\[-0.27ex]
\normalsize
 \vspace{4.27ex}{\Huge\BOZHO}	\\[4.27ex]
	\begin{tabular}{r@{\hspace{0.512em}}|@{\hspace{0.512em}}l}
 \vspace{0.27ex}	
&
 \vspace{0.27ex} 
	\end{tabular} \\[1.27ex]
 \vspace{0.27ex}	\\[0.27ex]
}

\begin{document}		

\renewcommand{\thepage}{\roman{page}}

\renewcommand{\thefootnote}{\fnsymbol{footnote}} 
\maketitle				
\renewcommand{\thefootnote}{\arabic{footnote}}   

\tableofcontents		


\begin{abstract}
This is a review with examples concerning the concepts of affine
(in particular, constant and linear) vector fields and fundamental vector
fields on a manifold. The affine, linear and constant vector fields on a
manifold are shown to be in a bijective correspondence with the fundamental
vector fields on it of respectively general affine, general linear and
translation groups (locally) represented on the manifold via the described in
this work left actions; in a case of the manifold
$\field^n=\field[R]^n,\field[C]^n$, the actions mentioned have the usual
meaning of affine, linear and translation transformations.
\end{abstract}

\renewcommand{\thepage}{\arabic{page}}


\section {Introduction}
\label{Introduction}

	This is a review with examples concerning the concepts of
affine (in particular, constant and linear) vector fields and fundamental
vector fields on a manifold.

	The linear vector fields are investigated in some detail in
section~\ref{Sect2}. On this ground, the constant and, more generally, affine
vector fields are briefly studied in section~\ref{Sect2-1}.
	Section~\ref{Sect3} is devoted to the flows of vector fields. As an
example, the flows of affine fields are found.
	Section~\ref{Sect3.3} is devoted to the canonical parameters and
invariants. As a major example, these functions for a constant vector field are
found, which can serve a guiding line for finding them for arbitrary affine
fields.
	Section~\ref{Sect4} deals with fundamental vector fields in the context
of affine vector fields, which turn to be such for particular representation of
the affine group of $\field[R]^n$. All of the results in the paper are
illustrated (or introduced at first) for the case of the manifold
$\field^n=\field[R]^n,\field[C]^n$.
\vspace{1ex}

	In this paper, the following notation will be used. By \field\ is
denoted the field of real or complex numbers, $\field=\field[R],\field[C]$.
A $C^1$ \field-manifold will be denoted by $M$, $(U,u)$ is a local chart of $M$
and $\{u^i\}$ is the coordinate system corresponding to it. The Latin indices
$i,j,k\dots$ run from~1 to~$\dim M$ and the Einstein's summation rule is
assumed. The space tangent to $M$ at $x\in M$ is $T_x(M)$. The set of vector
fields on $M$ is $\Fields{X}(M)$. If $f$ is a $C^1$ mapping between  $C^1$
manifolds, the induced by it tangent mapping is $T(f):=f_*$ and its restriction
to $T_x(M)$ is $T_x(f):=f_*|_x$. An arbitrary path with domain $J$,
$J\subseteq\field[R]$ being an arbitrary real interval, is denoted by $\gamma
\colon J\to M$.


\section{Linear vector fields} 	\label{Sect2}

	The importance of the linear vector fields comes from the facts that a
nonvanishing smooth vector field is locally linear and that their flows are
governed by simple equations (in suitable coordinates) which can be solved
explicitly. Usually, the linear vector fields are defined and considered in
$\field[R]^n$ equipped with standard Cartesian coordinates, as, e.g.,
in~\cite[p.~29]{Olver/LieGroups}, but below we are going to consider them on
arbitrary $C^1$ manifolds.

    \begin{Defn}    \label{Defn2.1}
Let $X$ be a vector field on a $C^1$ manifold $M$, $(U,u)$ be a local chart of
$M$ and $\{u^i\}$ be the local coordinates system associated with it. A
vector field $X\in\Fields{X}(M)$ with local representation
    \begin{equation}    \label{2.1}
X = X^i\frac{\pd }{\pd u^i}
    \end{equation}
is called \emph{linear  relative to $(U,u)$ or to $\{u^i\}$} if its local
components $X^i=X(u^i)$ relative to $\{u^i\}$ have a \field\ndash linear
dependence on the coordinate functions $u^1,\dots,u^{\dim M}$, \viz
    \begin{equation}    \label{2.2}
X^i = C_j^i u^j = \sum_{j=1}^{\dim M} C_j^i u^j
    \end{equation}
for some numbers $C_i^j\in\field$ forming the constant matrix $C:=[C_i^j]$.
    \end{Defn}

	The existence of linear (local) vector fields is evident: if we fix
$\{u^i\}$, then~\eref{2.1} and~\eref{2.2} define for a constant matrix $C$ a
vector field on $U$ which is linear relative to $\{u^i\}$. The interesting
problem is: which vector fields admit chart(s) with respect to which they are
linear?

    \begin{Prop}    \label{Prop2.1}
Let $X$ be a $C^1$ vector field on a $C^1$ manifold $M$, $p\in M$, and
$X_p\not=0$. Then there exists a chart $(U,u)$ such that $U\ni p$ and $X$ is
linear relative to $(U,u)$.
    \end{Prop}

    \begin{Proof}
According to~\cite[proposition~1.53]{Warner} or~\cite[p.~30]{Olver/LieGroups},
there is a chart $(V,v)$ of $M$ such that $V\ni p$ and
    \begin{equation}    \label{2.3}
X|_V = \frac{\pd }{\pd v^1} .
    \end{equation}
Consider a chart $(U,u)$ of the manifold $M$ with  $U\subseteq V$ and coordinate
functions
$u^i=f^i(v^2,\dots,v^{\dim M}) \e^{a^iv^1}$ (do not sum over $i$!) for some
suitable $C^1$ functions $f^i$, where $a^i\in\field$ are constants and
$a^1\not=0$.~%
\footnote{~%
Since we are interested only in the existence of the chart $(U,u)$, we do not
specify the functions $f^i$ and constants $a^i$. For instance, the choice
$a^1=1$, $f^1=a\in\field$, $a^k=0$ for $k\ge2$, and $f^k=v^k$ for $k\ge2$
results in the admissible coordinates $u^1=a\e^{v^1}$, $u^k=v^k$ for $k\ge2$
and $U=V$.%
}
According to~\eref{2.3}, we have
\[
X|_U
= \frac{\pd }{\pd v^1}\Big|_U
= \frac{\pd u^i}{\pd v^1}\Big|_U \frac{\pd }{\pd u^i}\Big|_U
= \sum_i a^i u^i \frac{\pd }{\pd u^i}  .
\]
Therefore $X$ is linear relative to the chart $(U,u)$ with local coordinates
$\{u^i\}$. In the particular case, the matrix $C$ is the diagonal matrix
$\diag(a^1,\dots,a^{\dim M})$.
    \end{Proof}

    \begin{Cor}    \label{Cor2.1}
If $X$ is a regular vector field on $M$, in a neighborhood of each point in $M$
there is a chart relative to which $X$ is linear.
    \end{Cor}

    \begin{Proof}
Since the regularity of $X$ means $X_p\not=0$ for all $p\in M$, the assertion
follows from proposition~\ref{Prop2.1}.
    \end{Proof}

    \begin{Exrc}    \label{Exrc2.1}
Prove that, if $X$ is linear relative to $(U,u)$, there is a chart $(V,v)$ with
$V\subseteq U$ in which~\eref{2.3} holds. (Hint: invert the proof of
proposition~\ref{Prop2.1}.)
    \end{Exrc}

    \begin{Exrc}    \label{Exrc2.2}
Show that, if $X$ is linear relative to $(U,u)$, then it is regular on some
open subset $V\subseteq U$, $X|_V\not=0$. (Hint: the equation $C_j^iu^j=0$
defines a $(\dim M-1)$\ndash dimensional submanifold of $U$ (hence of $M$).)
    \end{Exrc}

	If $X$ is linear relative to a coordinate system $\{u^i\}$, there exist
infinitely many other coordinates systems $\{v^i\}$ with respect to which it is
linear too. Indeed, if, for some $C_i^j,B_i^j\in\field$,
\[
X = (C_j^i u^j)\frac{\pd }{\pd u^i}
  = (B_j^i v^j)\frac{\pd }{\pd v^i}
\]
on the intersection of the domains of $\{u^i\}$ and $\{v^i\}$, then
    \begin{equation}    \label{2.4}
C_j^iu^j \frac{\pd v^k}{\pd u^i} = B_j^k v^j .
    \end{equation}
The problem is to be found $v^i$ and $B_i^j$ if $u^i$ and $C_i^j$ are given.
For instance, if $\{u^i\}\mapsto\{v^i\}$ is linear, \ie $v^i=a_j^iu^j$ with
$a_i^j\in\field$ and $\det[a_i^j]\not=0$, then the general solution
of~\eref{2.4} is $[B_i^j]=[a_i^j]\cdot C\cdot[a_i^j]^{-1}$. However,
equation~\eref{2.4} admits an infinite number of solutions
$(\{v^i\},\{B_j^i\})$ for which the change $\{u^i\}\mapsto\{v^i\}$ is
nonlinear.

	As a conclusion of the above results, we can say that the (local)
regularity of a $C^1$ vector field is equivalent to its linearity relative to
some local chart; such charts are not unique and their set is described via
equation~\eref{2.4}.

	Suppose a vector field $X\in\Fields{X}(M)$ is linear relative to a
chart $(U,u)$ and~\eref{2.1} and~\eref{2.2} hold. The vector fields
    \begin{equation}    \label{2.5}
E_i^j := u^j\frac{\pd }{\pd u^i} \in\Fields{X}(U)
    \end{equation}
form a basis on $U$ for the vector fields linear relative to $(U,u)$ in a sense
that any such field $X$ is on $U$ a linear combination of $E_i^j$ with
\emph{constant} coefficients, \viz
    \begin{equation}    \label{2.6}
X|_U = \sum_{i,j} C_j^i E_i^j .
    \end{equation}
The vector fields~\eref{2.5} are generators of a Lie algebra with respect to
the commutators as they are \field\ndash linearly independent and
    \begin{equation}    \label{2.7}
[ E_j^i , E_l^k ]
= ( \delta_j^k\delta_p^i\delta_l^q - \delta_l^i\delta_p^k\delta_j^q ) E_q^p
    \end{equation}
where $\delta_i^j$ are the Kronecker deltas.
Of course, the linear combinations $a_{jl}^{ik}E_k^l$, with
$a_{jl}^{ik}\in\field$ and $E_j^i\mapsto a_{jl}^{ik}E_k^l$ invertible, form
also a basis on $U$ for the vector fields linear with respect to $(U,u)$ and
generate a Lie algebra.

	Since the set $\bigl\{ \frac{\pd }{\pd u^i} \bigr\}$ of linearly
independent vector fields generates the module $\Fields{X}(M)$, the
$(\dim M)^2$ vector fields $\{E_j^i\}$ are linearly dependent for $\dim M\ge2$
and between them exist at least $(\dim M)^2-\dim M$ (independent) connections.


\section {Constant and affine vector fields}
\label{Sect2-1}

	Similar to definition~\ref{Defn2.1}, we define the constant and affine
vector field by
    \begin{Defn}    \label{Defn2.2}
 A vector field $X\in\Fields{X}(M)$ is termed \emph{constant} or \emph{affine}
(\emph{linear inhomogeneous}) relative to a chart $(U,u)$ of $M$, if in the
coordinate system $\{u^i\}$ associated with this chart, it has respectively the
representation
    \begin{align}    \label{2.8}
X &= B^i \frac{\pd }{\pd u^i}
\\		    \label{2.9}
X &= (C_j^i u^j + B^i) \frac{\pd }{\pd u^i}
    \end{align}
for some constant numbers $B^i,C_j^i\in\field$.
    \end{Defn}

	Obviously, the linear and constant vector fields are special case of
the affine ones for $B^i\equiv 0$ or $C_i^j\equiv 0$, respectively.

    \begin{Prop}[\normalfont \cf proposition~\ref{Prop2.1}\bfseries]
			\label{Prop2.2}
Let $X$ be a $C^1$ vector field on a $C^1$ manifold $M$, $p\in M$, and
$X_p\not=0$. Then there exists charts $(U,u)$ and $(W,w)$ such that $U\ni p$,
$W\ni p$ and $X$ is constant relative to $(U,u)$ and is affine with respect to
$(W,w)$.
    \end{Prop}

    \begin{Proof}
Let $(V,v)$ be a chart of $M$ such that $p\in V$ and~\eref{2.3} holds. Then, we
can put $(U,u)=(V,v)$ and define $(W,w)$ such that $W\subseteq V$ and assign
to it coordinates functions to be
 $w^i=f^i(v^2,\dots,v^{\dim M})\e^{a^iv^1} + B^i v^1$, where the notation of
the proof of proposition~\ref{Prop2.1} is used and $B^i\in\field$ are constant
numbers.
    \end{Proof}

	One can easily see (cf.~\eref{2.4}) that there exist infinitely many
charts relative to which a vector field is constant/affine if it is
constant/affine in a given chart; in particular, such coordinates systems for
affine (constant) vector fields can be obtained from a given one via linear
inhomogeneous transformations with \emph{constant} coefficients (resp.\ via
rescaling the coordinate  functions by constants).

	If $(U,u)$ is a chart of $M$, the two sets of vector fields
    \begin{align}    \label{2.10}
E_i &:= \frac{\pd }{\pd u^i}
\\		    \label{2.11}
E_i^j &:= u^j \frac{\pd }{\pd u^i} \quad
E_k = \frac{\pd }{\pd u^k}
    \end{align}
form a basis for the vector fields which are respectively constant or affine
relative to $(U,u)$ in a sense that any such field is a linear combination with
\emph{constant} coefficients of the above fields. The elements of the former
set generate an Abelian Lie algebra as $[E_i,E_j]=0$, while the ones of the
latter set generate a non-Abelian Lie algebra as their commutators are
    \begin{subequations}    \label{2.12}
    \begin{align}    \label{2.21a}
[E_i,E_j] &= 0
\\    \label{2.12b}
[E_i,E_k^j] &= \delta_i^j E_k
\\    \label{2.12c}
[ E_j^i , E_l^k ]
&=
( \delta_j^k\delta_p^i\delta_l^q - \delta_l^i\delta_p^k\delta_j^q ) E_q^p .
    \end{align}
    \end{subequations}


\section {Flows of vector fields}
\label{Sect3}

\subsection{The 2-dimensional case}
	\label{Subsect3.1}

	To begin with, consider a vector field $X$ on $\field[R]^2$ coordinated
by (generally said, non-Cartesian) coordinates $(u,v)$. We have the expansion
$X=x\frac{\pd }{\pd u}+y\frac{\pd }{\pd v}$ in which the functions $x$ and $y$
are the components of $X$ relative to the coordinate system $(u,v)$.~%
\footnote{~%
More precisely, $x$ and $y$ are the components of $X$ with respect to the
natural frame $\bigl\{\frac{\pd }{\pd u},\frac{\pd }{\pd v}\bigr\}$ induced by
the coordinates $(u,v)$.%
}
We have $X(u)=x$ and $X(v)=y$. If
$\beta_p \colon J\to \field[R]^2$, with $J$ being an open real interval,
$0\in J$ and $\beta_p(0)=p$ for some $p\in\field[R]$, is a path along which $X$
reduces to $\dot\beta$ (\ie $\beta$ is the integral path of $X$ through $p$ --
\emph{vide infra}), then the last equations reduce along $\beta$ to
    \begin{equation}    \label{3.1}
\frac{\od (u\circ\beta_p(t))}{\od t} = x\circ\beta_p (t)
\quad
\frac{\od (v\circ\beta_p(t))}{\od t} = y\circ\beta_p (t) ,
\qquad t\in J .
    \end{equation}
The solutions of these equations with respect to $\beta_p(t)$ such that
$\beta_p(0)=p$ for $p\in\field[R]^2$ define the integral paths of $X$. They
define a (local) 1\ndash parameter group $a$ of transformations of
the $\field[R]^2$ plane, termed also (local) flow of $X$, assigning to each
$t\in J$ the mapping
    \begin{equation}    \label{3.2}
    \begin{split}
a_t & \colon \field[R]^2\to\field[R]^2
\\
a_t & \colon \field[R]^2\ni p \mapsto a_t(p):=\beta_p(t)\in \field[R]^2
    \end{split}
    \end{equation}
and locally represented via the mapping
    \begin{equation}    \label{3.3}
(u_t,v_t)
:= (u,v)\circ a_t \circ (u,v)^{-1} \colon \field[R]^2\to\field[R]^2 .
    \end{equation}
So that, if a point $p\in\field[R]^2$ has coordinates  $(p_u,p_v)$ in $(u,v)$,
\ie $p=(u,v)^{-1}(p_u,p_v)$, then they change under the flow of $X$ according
to
\[
(u_t,v_t) \colon (p_u,p_v)\mapsto  (u,v)\circ a_t (p) =:(p_u(t),p_v(t)) .
\]

	Under the flow $a$ of a vector field $X$, the points in $\field[R]^2$
move along their orbits (the integral paths of $X$) and any function
$f \colon \field[R]^2\to \field[R]^2$ transforms into
    \begin{equation}    \label{3.4}
f_t:=f\circ a_t.
    \end{equation}
It is said that $f$ is dragged along $X$ (by the flow $a$) into $f_t$. If  $f$
is of class $C^1$, we have
    \begin{equation}    \label{3.5}
X(f)
= \frac{\od f_t}{\od t}\Big|_{t=0}
= \lim_{t\to 0} \frac{f_t-f}{t} .
    \end{equation}



\subsection{The general case}
	\label{Subsect3.2}

	Consider now the general case of a vector field $X\in\Fields{X}(M)$ on
a $C^1$ manifold $M$.

	In a local chart $(U,u)$ of $M$ with local coordinate system
$\{u^1,\dots,u^{\dim M}\}$ is valid the expansion
    \begin{equation}    \label{3.6}
X = X^i \frac{\pd }{\pd u^i}
  = \sum_{i=1}^{\dim M}  X^i \frac{\pd }{\pd u^i}
    \end{equation}
in which $X^i \colon U\to \field$ are the components of $X$ relative to the
(natural) frame $\Bigl\{\frac{\pd }{\pd u^i}  \Bigr\}$ on $U$. For a $C^1$
function $f\colon M\to\field$, we put
    \begin{equation}    \label{3.7}
f_i :=\frac{\pd f}{\pd u^i} \colon p\mapsto \frac{\pd f}{\pd u^i} \Big|_p
    := \frac{\pd (f\circ u^{-1})} {\pd r^i} \Big|_{u(p)}
    \end{equation}
with $\{r^1,\dots,r^{\dim M}\}$ being the standard Cartesian coordinate system
on $\field^{\dim M}$, that is
$r^i \colon \field^{\dim M}\ni(c_1,\dots,c_{\dim M}) \mapsto c_i$, so that
$u^i=r^i\circ u$ and $r^1=\id_{\field}$ for $\dim M=1$. Therefore
    \begin{align}    
			  \label{3.9}
X(u^i) = X^i .
    \end{align}

	Let $J$ be an open \field[R]-interval, $0\in J$, and $U$ be an open set
in $M$.  A \emph{local 1\ndash parameter group of local ($C^k$, if $M$ is of
class $C^k$) transformations}~\cite[ch.~I, \S~1]{K&N-1}
with domain $J\times U$ is a mapping
$a \colon J\times U\to M$,
$a \colon J\times U\ni (t,p) \mapsto a_t(p)\in M$, such that
	(i)~for any $t\in J$, the mapping
$a_t \colon U\ni p\mapsto a_t(p)\in a_t(U)$
is  $C^k$ diffeomorphism on $a_t(U)$ if $M$ is of class $C^k$ and
	(ii)~if $s,t,s+t\in J$ and $p,a_t(p)\in U$, then
    \begin{equation}    \label{3.9-1}
a_{s+t}(p) = a_s\circ a_t (p) \equiv a_s(a_t(p)) .
    \end{equation}

	If $U=M$ (and $J=\field[R]$), $a$ is called (global) 1\ndash parameter
group of transformations and $X\in\Fields{X}(M)$ is said to be complete in
that case.

	A local 1-parameter group $a$ of transformations induces on $U$ a
vector field $X$ via the equation $X_p=\dot{\beta}_p(0)$, where $p\in U$ and
the path $\beta_p \colon J\to M$ is given by
    \begin{equation}    \label{3.9-2}
\beta_p(t):=a_t(p) ,
    \end{equation}
\ie $\beta_p$ is the orbit of $p$ and, at the same time, the integral path of
$X$ trough $p$.

	The inverse is also true (see~\cite[ch.~I, \S~1,
proposition~1.5]{K&N-1}) or~\cite[theorem~1.48]{Warner}). If
$X\in\Fields{X}(M)$, for every $p\in M$ there exist an interval $J_p\ni 0$,
open set $U_p$, and local 1\ndash parameter group $a \colon J_p\times U_p\to M$
of local transformations which induces the given vector field $X$ (on $U_p$).
This local group is called \emph{(local) flow} of $X$. In local coordinates it
can be defined as follows.

	Suppose $(U,u)$ is a chart of $M$ and $\{u^i=r^i\circ u\}$ is the
corresponding local coordinate system. Recall that $\beta \colon J\to M$ is an
integral path of $X$ iff $X\circ\beta=\dot\beta$.
Hence, if $X=X^i\frac{\pd }{\pd u^i}$, we get
\(
X_{\beta(t)} (u^i)
= (\dot{\beta}(t))(u^i)
= \frac{\od (u^i\circ\beta(t))}{\od t} .
\)
Combining the last result with~\eref{3.9}, we find the equation of the
integral paths of $X$ as
    \begin{equation}    \label{3.10}
\frac{\od (u^i\circ\beta(t))}{\od t} = X^i\circ\beta(t) .
    \end{equation}

	Let $\beta_p \colon J\to M$ be the integral path of $X$ through a point
$p\in M$, \ie $\beta_p(0)=p$ and $\beta_p$ satisfies~\eref{3.10} with
$\beta_p$ for $\beta$. Then the flow of $X$ is given by $a_t(p)=\beta_p(t)$ and
in $(U,u)$ is represented by
    \begin{equation}    \label{3.10-1}
u_t := u\circ a_t\circ u^{-1} \colon \field^{\dim M}\to \field^{\dim M}
    \end{equation}
and satisfies the system of differential equations

    \begin{subequations}    \label{3.11}
    \begin{gather}    \label{3.11a}
\frac{\od u_t^i(\lambda)}{\od t}
=
X^i\circ u^{-1}\circ u_t(\lambda)
=
X^i\circ u^{-1}\circ (u_t^1,\dots,u_t^{\dim M})(\lambda)
\intertext{for $i=1,\dots,\dim M$ and $\lambda\in u(U)$, due to~\eref{3.10},
and the initial condition}
		    \label{3.11b}
(u_t^1,\dots,u_t^{\dim M})|_{t=0} = \id_{\field^{\dim M}} .
    \end{gather}
    \end{subequations}

	Equations~\eref{3.11a} simply mean that~\eref{3.9} are the local
equations of the flow of $X$ if $X^i=f^i(u^1,\dots,u^{\dim M})$ are known
expressions of the local coordinate functions $u^1,\dots,u^{\dim M}$. Indeed,
if this is the case, equations~\eref{3.9} along the integral paths of $X$
reduce to (see~\eref{3.10} and~\eref{3.4})
    \begin{equation}
	\tag{\ref{3.11a}$^\prime$}	\label{3.11a'}
\frac{\od u_t^i}{\od t} = f^i(u_t^1,\dots,u_t^{\dim M}) ,
    \end{equation}
which is an equivalent form of~\eref{3.11a}

	Under the flow $a$ of $X\in\Fields{X}(M)$, a function
$f \colon M\to \field$ is dragged according to~\eref{3.4} and, if $f$ is of
class $C^1$, then its directional derivative with respect to $X$ is given
by~\eref{3.5}.

\subsection{Affine vector fields}
	\label{Subsect3.4}

	Let us now find the flow of an affine vector field with local
representation~\eref{2.9}. For the purpose, we shall use the following matrix
notation
    \begin{equation}    \label{3.31}
    \begin{split}
U &:= (u^1,\dots,u^{\dim M})^\top \quad
\frac{\pd }{\pd U} :=
		\Bigl( \frac{\pd }{\pd u^1},\dots,\frac{\pd }{\pd u^{\dim M}}
		\Bigr)
\\
C &:= [C_j^i] \quad
B := (B^1,\dots,B^{\dim M})^\top
    \end{split}
    \end{equation}
in which the affine field~\eref{2.9} has the form~%
\footnote{~%
At this point, we use the rule that summation excludes differentiation as a
result of which the r.h.s.\ of~\eref{3.32} is simply a shortcut for the r.h.s.\
of~\eref{2.9}.%
}
    \begin{equation}    \label{3.32}
X = \frac{\pd }{\pd U}\cdot (C\cdot U+B)
    \end{equation}
and hence
    \begin{equation}    \label{3.33}
X(U) = C\cdot U+B .
    \end{equation}

	In accordance with~\eref{3.10-1}, the flow $a$ of $X$ is represented in
the coordinates $\{u^i\}$ by the matrix $U_t=(u_t^1,\dots u_t^{\dim M})^\top$
with $u_t^i=u^i\circ a\circ((u^1)^{-1},\dots,(u^{\dim M})^{-1})$, \ie
 $U_t=U\circ a_t\circ((u^1)^{-1},\dots,(u^{\dim M})^{-1})$.
According to~\eref{3.14} and~\eref{3.33}, $U_t$ is the solution of the
following matrix initial\ndash valued problem:
    \begin{equation}    \label{3.34}
\frac{\od U_t}{\od t} = C\cdot U_t + B
\quad U_t|_{t=0} = U .
    \end{equation}
Therefore the explicit form of $U_t$ is
    \begin{subequations}    \label{3.35}
    \begin{alignat}{3}    \label{3.35a}
U_t &= tB+U 	&&\qquad \text{for } C=0 &&\text{ (constant field)}
\\		    \label{3.35b}
U_t &= \e^{tC}U	&&\qquad \text{for } B=0 &&\text{ (linear field)}
\\		    \label{3.35c}
U_t &= \e^{tC}(U-\tilde{U})+\tilde{U}
		&&\qquad \text{for } C\not=0 \text{ and } B\not=0
					 &&\text{ (affine field)} ,
    \end{alignat}
    \end{subequations}
where the constant matrix $\tilde{U}$ is a solution of the equation
    \begin{equation}    \label{3.36}
C\cdot\tilde{U} + B = 0
\qquad\text{(with $C\not=0$ and $B\not=0$)}.
    \end{equation}
(Note, the r.h.s.\ of~\eref{3.35c} is independent of the particular choice of
$\tilde{U}$ if equation~\eref{3.36} has more than one solution with respect to
$\tilde{U}$; if $C$ is nondegenerate, then $\tilde{U}=-C^{-1}\cdot B$.) Thus the
flow of a constant vector field is governed by a linear law while the one of a
linear (affine) vector field  is governed by an exponential law (combined with
shift by $\tilde{U}$).



\section{Invariants and canonical parameters}
	\label{Sect3.3}

    \begin{Defn}    \label{Defn3.1}
	A $C^1$ function $I \colon M\to \field$ is called an \emph{invariant}
of $X\in\Fields{X}(M)$ if it is constant along the integral paths of $X$,
$I_t=I=\const$ for all $t$ or $X(I)=0$.
	A $C^1$ function $S \colon M\to \field$ is termed a \emph{canonical
parameter} of $X$ if $X(S)=1$.
    \end{Defn}

	The following result is completely obvious but worth recording.

    \begin{Prop}    \label{Prop3.1}
The difference of two canonical parameters of $X$ is an invariant of $X$ and
the sum of a canonical parameter and invariant of $X$ is a canonical parameter
of $X$. Consequently, a canonical parameter is defined up to an invariant.
    \end{Prop}

	Any constant function $M\to\{c\}$ for a given $c\in\field$ is an
invariant of all vector fields on $M$. However, the existence of non\ndash
trivial invariants as well as of canonical parameters is not evident.
Generally they exist only locally as stated in the following result.

    \begin{Prop}    \label{Prop3.2}
Given a point $p\in M$ and  $C^1$ vector field on a $C^3$ real manifold $M$ (or
$C^3$ complex manifold considered as a real one of real dimension
$\dimR M=2\dimC M$) such that $X_p\not=0$. Then there is an open set
$V\subseteq M$ containing $p$, $V\ni p$, on which exist a canonical parameter
$S$ and non\ndash constant invariant $I$ of $X$.
    \end{Prop}

    \begin{Proof}
According to~\cite[proposition~1.53]{Warner}, there is a local chart $(V,v)$
with coordinate functions $v^i$ such that $V\ni p$ and
    \begin{equation}    \label{3.14}
X|_V = \frac{\pd }{\pd v^1}  .
    \end{equation}
Defining
    \begin{subequations}    \label{3.15}
    \begin{align}    \label{3.15a}
S = v^1 + F(v^2,\dots,v^{\dim M})
\\		    \label{3.15b}
I = G(v^2,\dots,v^{\dim M}) ,
    \end{align}
    \end{subequations}
where $F$ and $G$ map a $(\dim M-1)$-tuple of $C^1$ functions on $V$ into a
$C^1$ function on $V$, from~\eref{3.14}, we see that $S$ and $I$ are
respectively a canonical parameter and an invariant of $X$ on $V$, \ie
$X|_V(S)=1$ and $X|_V(I)=0$.
    \end{Proof}

    \begin{Cor}    \label{Cor3.1}
Under the hypotheses of proposition~\ref{Prop3.2} and the notation introduce
in its proof, all local canonical parameters and invariants of a vector field
are given via~\eref{3.15a} and~\eref{3.15b}, respectively.
    \end{Cor}

    \begin{Proof}
Use proposition~\ref{Prop3.1} and the proof of proposition~\ref{Prop3.2}.
    \end{Proof}

    \begin{Exrc}    \label{Exrc3.1}
If $\dim M=1$, prove that the only (\field\ndash valued) invariants of a vector
field, with separable (by open sets) points at with it is irregular, if any,
are the constant functions $M\to\{c\}$ for some $c\in\field$.
    \end{Exrc}

	In a neighborhood of a point $p$ at which $X_p\not=0$, the proof of
proposition~\ref{Prop3.2} provides the canonical parameter $v^1$ and, if
$\dim M\ge 2$, the $n-1$ invariants $v^2,\dots,v^{\dim M}$. The set of these
functions $\{v^i\}$ is a coordinate system on $V$ in which~\eref{3.14}
holds. The converse of that observation reads

    \begin{Prop}    \label{Prop3.3}
Let $X\in\Fields{X}(M)$, $S$ be a canonical parameter of $X$, and
$I^2\dots,I^{\dim M}$ be invariants of $X$ on $U\subseteq M$. If the set
$\{S,I^2\dots,I^{\dim M}\}$  is a coordinate system on $U$, \ie for some chart
$(U,u)$ of $M$, then in it
    \begin{equation}    \label{3.16}
X|_U = \frac{\pd }{\pd S}
    \end{equation}
and the flow $a$ of $X$ in it is represented via the mapping
$u\mapsto u_t=(r^1+t,r^2,\dots,r^{\dim M})$ (see~\eref{3.10-1})  or
    \begin{equation}    \label{3.17}
u_t \colon (s,\bs{i})\mapsto (s_t,\bs{i}_t) = (s+t,\bs{i})
    \end{equation}
for all $s\in\field$ and $\bs{i}\in\field^{\dim M -1}$.
    \end{Prop}

    \begin{Proof}
Let $\{u^i\}$ be  arbitrary coordinate system on $U$ and
$X|_U=X^i\frac{\pd }{\pd u^i}$. Making the change
$\{u^i\}\mapsto \{S,I^1\dots,I^{\dim M}\}$, we get
\[
X|_U
= X^i\frac{\pd }{\pd u^i}
=
X^i\Bigl\{
\frac{\pd S}{\pd u^i} \frac{\pd }{\pd S}
+\sum_{k=2}^{\dim M} \frac{\pd I^k}{\pd u^i} \frac{\pd }{\pd I^k}
 \Bigr\}
=
X(S)\frac{\pd }{\pd S} +\sum_{k=2}^{\dim M} X(I^k) \frac{\pd }{\pd I^k}
=
\frac{\pd }{\pd S} ,
\]
where $X(f) = f_i X^i$, for a $C^1$ function $f$, and definition~\ref{Defn3.1}
were applied. To prove~\eref{3.17}, we notice that, by virtue of~\eref{3.16},
the equations~\eref{3.11a} in $\{S,I^1\dots,I^{\dim M}\}$ read
\[ \frac{\od
u_t^1(\lambda)}{\od t} = 1 \quad
\frac{\od u_t^k(\lambda)}{\od t} = 0 \text{ for } k\ge 2
\]
and their solution $u_t=(u_t^1,\dots,u_t^{\dim M})$, under the
condition~\eref{3.11b}, is $u_t^1=r^1+t$ and $u_t^k=r^k$ for $k>1$, where
$\{r^i\}$ is the standard coordinate system on $\field^{\dim M}$ and $t$ is
considered as the constant function $\field^{\dim M}\to t$.
    \end{Proof}

\begin{Exmp}\label{Exmp3.2}
	Consider $\field[R]^2$ coordinated by the standard Cartesian coordinates
$(u,v)=(r^1,r^2)$ and the affine vector field
$X=\alpha\frac{\pd }{\pd u}+(2\beta u + \gamma)\frac{\pd }{\pd v}$ for some
$\alpha,\beta,\gamma\in\field[R]$ with $\alpha\not=0$. The flow $a$ of $X$ is
locally represented by (see~\eref{3.10-1})
\(
(u_t,v_t) := (u,v)\circ a_t \circ(u,v)^{-1}  \colon \field[R]^2\to \field[R]^2
\)
and is the solution of the initial-value problem (see~\eref{3.11})
    \begin{align*}
 \frac{\od u_t}{\od t}
& = \alpha \quad \frac{\od v_t}{\od t}
  = (2\beta u +\gamma)\circ (u,v)^{-1}\circ(u_t,v_t)
  = 2\beta u_t + \gamma ,
\\
 (u_0,v_0) &= (r^1,r^2) = \id_{\field^{2}} ,
    \end{align*}
so that $u_t=u+\alpha t$ and $v_t=v+(2\beta u+\gamma)t+\alpha\beta t^2$, \ie
 \[
 (u_t,v_t) \colon \field[R]\ni(b,c) \mapsto
		(b_t,c_t)=(b+\alpha t,c+(2\beta b+\gamma)t+\alpha\beta t^2).
 \]
From here and~\eref{3.9-2}, we see that the local coordinates of the point
$\beta_p(t)$ of the integral path of $X$ passing through $p=(u,v)^{-1}(b,c)$
are
	\begin{align*}
      (u,v)(\beta_p(t))
 & = (b_t,c_t)=(b+\alpha t,c+(2\beta b+\gamma)t+\alpha\beta t^2)
	\end{align*}

	The particular vector field $X$ has a (global) canonical parameter
$S=\frac{1}{\alpha}u$ and invariant $I=\alpha v-\beta u^2-\gamma$, $X(S)=1$ and
$X(I)=0$. This canonical parameter and invariant (as well as all of them) can
be found in a way similar to the one described below in example~\ref{Exmp3.12}.

	Applying corollary~\ref{Cor3.1}, we can assert that all canonical
parameters of $X$ are $S+F(I)$ and all its invariants are $G(I)$, where $F$ and
$G$ map $C^1$ functions on $\field[R]^2$ into $C^1$ functions on $\field[R]^2$

	The Jacobian of the coordinate change $(u,v)\mapsto(S,I)$ is
\(
    \bigl|\begin{smallmatrix}
1/ \alpha  & 0 \\
-2\beta u   & \alpha
    \end{smallmatrix}\bigr|
=1\not=0 .
\)
Therefore the pair $(S,I)$ induces the natural frame
$\bigl\{\frac{\pd }{\pd S},\frac{\pd }{\pd I}  \bigr\}$
on the set $U=\field[R]^2$ and, by proposition~\ref{Prop3.3}, on $U$ in
$(S,I)$ we have
\[
X|_U = \frac{\pd }{\pd S}
\qquad
(S_t,I_t) \colon U\ni(b,c )\mapsto (b_t,c_t) = (b+t,c) .
\]
\EndExmp
\end{Exmp}

    \begin{Exmp}    \label{Exmp3.12}
Consider a constant vector field $X$  with nonvanishing coefficients,
\[
X= B^i\frac{\pd }{\pd u^i} \in \Fields{X}(M)
\qquad B^i\not=0 \text{ for all } i .
\]
The $C^1$ functions $S$ and $I$ are respectively a canonical parameter and an
invariant of $X$ iff they are solutions of the differential equations
\[
B^i\frac{\pd S}{\pd u^i} = 1 , \quad
B^i\frac{\pd I}{\pd u^i} = 0
\]
According to the general theory of differential equations of this
kind~\cite[pp.~733--735]{Matveev1974}, these equations, relative to $S$ and
$I$, are equivalent to respectively the systems
    \begin{align*}
\frac{\od u^1}{B^1} & = \dots =\frac{\od u^{\dim M}}{B^{\dim M}}
=\frac{\od S}{1}
\\
\frac{\od u^1}{B^1} & = \dots =\frac{\od u^{\dim M}}{ B^{\dim M}}
\quad \od I =0
    \end{align*}
and, consequently, the general form of the functions $S$ and $I$ is determined
via the equations
 $\Phi(\varphi_1,\dots,\varphi_{\dim M}) = 0$ and
 $\Psi(\psi_1,\dots,\psi_{\dim M}) = 0$,
where $\Phi$ and $\Psi$ are arbitrary $C^1$ functions and
    \begin{alignat*}{2}
\varphi_1 & = S - \frac{u^1}{B^1} &&\quad
\varphi_i = u^i - u^1\frac{B^i}{B^1} \text{ for } i\ge2
\\
\psi_1 & = I &&\quad
\psi_i = u^i - u^1\frac{B^i}{B^1} \text{ for } i\ge2
    \end{alignat*}
are  $n$ independent integrals of the above systems. Admitting that the last
equations can be solved with respect to $S$ and $I$, we get
    \begin{align*}
S &=
\frac{u^1}{B^1}
+ F\Bigl(
u^2- u^1\frac{B^2}{B^1}, \dots,
u^{\dim M} - u^1\frac{B^{\dim M}}{B^1}
\Bigr)
\\
I &=
   G\Bigl(
u^2- u^1\frac{B^2}{B^1}, \dots,
u^{\dim M} - u^1\frac{B^{\dim M}}{B^1}
\Bigr)
    \end{align*}
for some $C^1$ functions $F$ and $G$, which agrees with corollary~\ref{Cor3.1}.
\EndExmp
    \end{Exmp}

    \begin{Rem}    \label{Rem3.1}
If some of the $B$'s in~\eref{2.8} vanish, the procedure described in
example~\ref{Exmp3.12} remains valid for the remaining non-vanishing $B$'s and
to the obtained in this way expressions for $S$ and $I$ can be added arbitrary
functions of the $u$'s for which the similarly indexed $B$'s vanish.
    \end{Rem}

	In a way similar to the one considered in example~\ref{Exmp3.12}, one
can investigate the problems of finding the canonical parameters and invariants
of linear and affine vector fields. The only difference  from the case of a
constant vector field is in the more difficult differential equations that
should be solved.



\section {Fundamental vector fields}
\label{Sect4}

	A left (right) action of a Lie group $G$, with multiplication
$G\times G\ni(a,b)\mapsto ab\in G$ and identity element $e$, on a manifold $M$
is a mapping
 $\lambda \colon G\times M\to M$ ($\rho \colon M\times G\to M$)
such that
  $\lambda_{ab}=\lambda_a\circ\lambda_b$ and $\lambda_e=\id_M$
($\rho^{ab}=\rho^b\circ\rho^a$ and $\rho^e=\id_M$) for $a,b\in G$,
where the partial mappings
 $\lambda_a,\rho^a \colon M\to M$ are defined by
$\lambda_a(x):=\lambda(a,x)$ and $\rho^a(x):=\rho(x,a)$ for all $x\in M$.
	Below we shall need also the partial mappings
 $\lambda^x,\rho_x \colon G\to G$ defined by
$\lambda^x(a):=\lambda(a,x)$ and $\rho_x(a):=\rho(x,a)$ for all $x\in M$ and
$a\in G$.
	A fundamental vector field on $M$ is a vector field on $M$ which is
obtained from a vector in the space $T_e(G)$ tangent to $G$ at the identity
element $e\in G$ via the tangent mappings of $\lambda^x$ or $\rho_x$.
Precisely, we have following definition.

    \begin{Defn}    \label{Defn4.1}
Let $\lambda$ and $\rho$ be respectively left and right actions of a $C^1$ Lie
$G$ on $M$. For $X\in T_e(G)$,
the \emph{left and right fundamental vector fields}
$\xi_X^l\in\Fields{X}(M)$ and $\xi_X^r\in\Fields{X}(M)$ on $M$ (associated
with $X$ and the given actions of $G$) are defined respectively by
    \begin{align*}
\xi_X^l \colon x\mapsto
\xi_X^l(x) &:= (T_e(\lambda^x))(X) = (T_{(e,x)}(\lambda))(X,0_x)
\in T_{\lambda(e,x)}(M) = T_{x}(M)
\\
\xi_X^r \colon x\mapsto
\xi_X^r(x) &:= (T_e(\rho_x))(X) = (T_{(x,e)}(\rho))(0_x,X)
\in T_{\rho(x,e)}(M) = T_{x}(M)
    \end{align*}
where $0_x\in T_x(M)$ is the zero vector in the space $T_x(M)$ tangent to $M$
at $X\in M$.
    \end{Defn}

	The above definitions can be localized in an evident way if we replace
in them $M$ by an arbitrary open subset $V\subseteq M$.

	Obviously, the vectors $\xi_X^l(x)$ and $\xi_X^r(x)$ are tangent at $x$
to the orbits $\lambda^x(G)$ and $\rho_x(G)$, respectively, of $G$ through $x$.

	More information on fundamental vector fields can be found, for
instance, in~\cite[pp.~46--47]{KMS-1993}, \cite[pp.~283--284]{Poor},
\cite[pp.~121--124]{Greub&et_al.-2} and~\cite{K&N-1}. It should be noted that
sometimes, \eg in~\cite{Chebotarev-1940}, \cite{Lie/Groups}
and~\cite[see especially pp.~127--133]{Cantwell}, the term `operator of a group'
is used instead of the modern one `fundamental vector field' associated to a
group (action on a manifold).

	Since the sets of left and right invariant vector fields,~%
\footnote{~%
If $L_a \colon G\ni b\mapsto ab$ and $R_a \colon G\ni b\to ba$ are the left and
right translations on $G$ by an element $a\in G$, a vector field
$X\in\Fields{X}(G)$ is left or right invariant if
$(L_a)_*(X_b)=X_{L_ab}=X_{ab}$  or
$(R_a)_*(X_b)=X_{R_ab}=X_{ba}$, respectively.%
}
on $G$ are isomorphic as vector spaces to $T_e(G)$~%
\footnote{~\label{isomorphism}%
If $X\in\Fields{X}(G)$ is left or right invariant, the isomorphism
$I \colon \{X\}\to T_e(G)$ is given by $I \colon X\mapsto X_e$, sending $X$ to
its value at the identity $e$, as
 $X_a=X_{L_e(a)}=(L_a)_*(X_e)$ or
 $X_a=X_{R_e(a)}=(R_a)_*(X_e)$, respectively.%
}
the left (resp.\ right) fundamental vector fields are images of the left or
right invariant vector fields on $G$ via the left (resp.\ right) action of the
group;~%
\footnote{~%
Indeed, if $Y$ is a left or right invariant vector field on $G$, the left
action $\lambda \colon G\times M\to M$ sends it into a vector field
$Y^l\in\Fields{X}(M)$ such that
$Y^l \colon X\mapsto Y^l(x):=Y^l(\lambda(e,x))=(T_e(\lambda^x))(Y|_e)$,
which depends only on $Y$ at the identity $e$ of $G$.
	Similarly, a right action $\rho \colon M\times G\to M$ sends $Y$ to
$Y^r\in\Fields{X}(M)$ such that
 $Y^r \colon x\mapsto Y^r(x)=(T_e(\rho_x))(Y|_e)$.
}
note that a \emph{left/right} fundamental vector fields are connected
with the \emph{left/right action of the group} and not with the left or right
invariant vector field from which one has started.~%
\footnote{~%
The reason for that being the isomorphism $I$ described in
footnote~\ref{isomorphism} above:
 $\xi_{I(X)}^l(x)=T_e(\lambda^x)(X_e)$ or
 $\xi_{I(X)}^r(x)=T_e(\rho_x)(X_e)$
regardless is $X\in\Fields{X}(M)$ left or right invariant.%
}
Consequently the fundamental vector fields are images of the Lie algebra $\g$
of $G$ via the group actions as $\g$ is identified with the set of left
invariant vector fields on $G$~\cite{K&N-1}, or with the set of right invariant
vector fields on $G$~\cite[p.~42, definition~1.44]{Olver/LieGroups}, or with
the space $T_e(G)$ tangent to $G$ at the identity element $e$~\cite{KMS-1993}.

	If $\{y^\mu:\mu=1,\dots,\dim G\}$ are local coordinates on an open
subset $U\subseteq G$ and $\{u^k:k=1,\dots,\dim M\}$ are coordinates on $M$
with $\lambda^x(U)$ in their domain, then~\cite[sec.~1.23(a)]{Warner}
    \begin{equation}    \label{4.0-1}
(T_a(\lambda^x))\Bigl( \frac{\pd }{\pd y^\mu}\Big|_a \Bigr)
=
\frac{\pd (u^k\circ\lambda^x)}{\pd y^\mu} \Big|_a
\frac{\pd }{\pd u^k} \Big|_{\lambda^x(a)} \qquad a\in G ;
    \end{equation}
we have a similar equation with $\rho_x$ for $\lambda^x$ in a case of a right
action $\rho$ when $\rho_x(U)$ is in the domain of $\{u^k\}$. If $e\in U$ and
$X=X^\mu\frac{\pd }{\pd y^\mu}\big|_e$, then the last equation immediately
implies
    \begin{subequations}    \label{4.0-2}
    \begin{align}    \label{4.0-2a}
\xi_X^l(x)
= X^\mu \frac{\pd (u^k\circ\lambda^x)}{\pd y^\mu} \Big|_e
	\frac{\pd }{\pd u^k}\Big|_{x}
\\		    \label{4.0-2b}
\xi_X^r(x)
= X^\mu \frac{\pd (u^k\circ\rho_x)}{\pd y^\mu} \Big|_e
	\frac{\pd }{\pd u^k}\Big|_{x} .
    \end{align}
    \end{subequations}


\subsection{Fundamental vector fields on \protect{$\field[R]^n$} induced
	    by \protect{$\GL(n,\field[R])$}}
	\label{Subsect4.1}

	Let the general linear group $\GL(n,\field[R])$, consisting of all
regular $n\times n$, $n\in\field[N]$, matrices with matrix multiplication as a
group multiplication and the identity $n\times n$ matrix $\openone$ as an
identity element, be represented on $\field[R]^n$
via a left action $\lambda \colon \GL(n,\field[R])\times\field[R]^n\to
\field[R]^n$ given by left multiplication, viz., if
$a=[a_j^i]\in\GL(n,\field[R])$ and $x=(x^1,\dots,x^n)^\top\in\field[R]^n$,
where $\top$ means matrix transposition, then
$\lambda_a(x)=\lambda^x(a)=\lambda(a,x)=a\cdot x=(a_j^1x^j,\dots,a_j^nx^j)^\top$.

    \begin{Prop}    \label{Prop4.1}
The vector field
    \begin{equation}    \label{4.4}
\xi_X \colon x\mapsto \xi_X(x)
:= T_{\openone}(\lambda^x)(X)
= \sum_{i,j} X_j^i u^j(x) \frac{\pd }{\pd u^i}\Big|_x
\in T_x(\field[R]^n)
    \end{equation}
is a (left) fundamental vector field on $\field[R]^n$ corresponding to
$X\in T_{\openone}(\GL(n,\field[R]))$ and the represention $\lambda$ of
$\GL(n,\field[R])$ on $\field[R]^n$ by left multiplication. Here $X_i^j$ are
the components of $X$ in the coordinate system defined a little below in the
proof of this assertion.
    \end{Prop}

    \begin{Proof}
	Consider the tangent mapping
\(
T_a(\lambda^x) \colon T_a(G)\to
T_{\lambda^x(a)}(\field[R]^n) = T_{a\cdot x}(\field[R]^n) .
\)
In the global coordinate systems $\{y_j^i\}$ on $\GL(n,\field[R])$ and
$\{u^i\}$ on $\field[R]^n$, given by $y_j^i(a):=a_j^i$ and $u^i(x):=x^i$, it is
represented by a Jacobi matrix with elements
$\frac{\pd (u^k\circ\lambda^x)}{\pd y_j^i}$ such
that
\[
(T_a(\lambda^x))\Bigl( \frac{\pd }{\pd y_j^i}\Big|_a \Bigr)
=
\frac{\pd (u^k\circ\lambda^x)}{\pd y_j^i} \Big|_a
\frac{\pd }{\pd u^k} \Big|_{\lambda^x(a)}
\]
due to ~\eref{4.0-1}. Since
    \begin{equation}    \label{4.0}
u^k\circ\lambda^x \colon a\mapsto
u^k(\lambda^x(a)) = u^k(a\cdot x) =a_j^kx^j = u^j(x) y_j^k(a)
    \end{equation}
the elements of the Jacobi matrix are $u^j(x)\delta_i^k$, so that
    \begin{equation}    \label{4.1}
(T_a(\lambda^x))\Bigl( \frac{\pd }{\pd y_j^i}\Big|_a \Bigr)
=
u^j(x) \frac{\pd }{\pd u^i}\Big|_{\lambda^x(a)}
= x^j \frac{\pd }{\pd u^i} \Big|_{a\cdot x}
    \end{equation}
and, if
\(
Y_a
= (Y_a)_j^i\frac{\pd }{\pd y_j^i}\big|_a ,
\)
then
    \begin{equation}    \label{4.2}
(T_a(\lambda^x))(Y_a)
= (Y_a)_j^i u^j(x) \frac{\pd }{\pd u^i}\Big|_{\lambda^x(a)}
= (Y_a)_j^i x^j \frac{\pd }{\pd u^i}\Big|_{a\cdot x} .
    \end{equation}
The particular settings $a=\openone=[\delta_j^i]$ and
$X=Y_{\openone}\in T_{\openone}(\GL(n,\field[R]))$ in the last equation reduce
it to
    \begin{equation}    \label{4.3}
(T_{\openone}(\lambda^x))(X) = X_j^i u^j(x) \frac{\pd }{\pd u^i}\Big|_x
    \end{equation}
from where~\eref{4.4} follows.
    \end{Proof}

	Evidently, the fundamental vector field
    \begin{equation}    \label{4.5}
\xi_X = X_j^i u^j \frac{\pd }{\pd u^i} \in \Fields{X}(\field[R]^n)
    \end{equation}
is a linear vector field (in the frame/coordinates used above) to which
corresponds the matrix $C=[X_j^i]$. In particular, to the vectors in
$T_{\openone}(\GL(n,\field[R]))$ with components $X_l^k=\delta_l^i\delta_j^k$
in $\{y_j^i\}$ correspond the fundamental vector fields
    \begin{equation}    \label{4.6}
 E_j^i = u^i \frac{\pd }{\pd u^j} \in \Fields{X}(\field[R]^n)
    \end{equation}
which form a basis for the set of linear vector fields in the sense that any
linear vector field is a linear combination with \emph{constant} coefficients
of these vector fields. These vector fields are generally linearly dependent
and between them exist at least $n^2-n$ (independent) connections.

    \begin{Cor}    \label{Cor4.1}
All linear (relative to Cartesian coordinates) vector fields on $\field[R]^n$
are fundamental vector fields of $\GL(n,\field[R])$ represented on
$\field[R]^n$ via left multiplication and \emph{vice versa}.
    \end{Cor}

    \begin{Proof}
If
 $c_j^i u^j \frac{\pd }{\pd u^j}\in \Fields{X}(\field[R])^n$ is a linear vector
field in the Cartesian coordinate system $\{u^i\}$, then, by
equation~\eref{4.5}, it is the fundamental vector field corresponding to
the vector
 $X=c_j^i \frac{\pd }{\pd y_j^i}\in T_{\openone}(\GL(n,\field[R]))$.
The converse was proved above.
    \end{Proof}

    \begin{Rem}    \label{Rem4.1}
A vector field on $\field[R]^n$ which is linear relative to non\ndash Cartesian
coordinates need not to be a fundamental vector field for $\GL(n,\field[R])$.
For instance, if $n=1$ and $u$ is the standard Cartesian coordinate on
$\field[R]$, the field $X=\frac{u}{u+1}\frac{\od }{\od u}$ is non\ndash linear
in $\{u\}$ and, consequently, is non-fundamental for $\GL(1,\field[R])$, but in
a non\ndash Cartesian coordinate system with the coordinate function $t=u\e^u$
it has the representation $X=t\frac{\od }{\od t}$ and hence it is linear in
$\{t\}$.
	The choice of the coordinates $\{y_i^j\}$ on $\GL(n,\field[R])$ is not
so important. According to~\eref{4.0}, their change results in replacing $x^j$
in~\eref{4.1} (or $u^j$ in or after it) with $A_i^jx^i$ (resp.\ $A_i^ju^i$)
where $A_i^j$ are constants depending only on the matrix $a$. Similar is the
result if instead of the standard Cartesian coordinates $\{u^i\}$ on
$\field[R]^n$ one uses any set of non\ndash standard Cartesian coordinates on
$\field[R]^n$.

	Similar remarks hold true with respect to Corollaries~\ref{Cor4.2}
and~\ref{Cor4.3} below.
    \end{Rem}



\subsection{Fundamental vector fields on \protect{$\field[R]^n$} induced
	    by \protect{$\T_n$}}
	\label{Subsect4.2}

	As a set, the translation group $\T_n$ on $\field[R]^n$ coincides with
$\field[R]^n$, $\T_n=\field[R]^n$, with addition as a group multiplication and
the zero vector as identity element;
hence $\T_n$ is an Abelian group. Its left and right actions
 $\lambda\colon \T_n\times\field[R]^n\to \field[R]^n$ and
 $\rho   \colon \field[R]^n\times \T_n\to \field[R]^n$,
respectively, on $\field[R]^n$ are defined by
 $\lambda(t,x)=\rho(x,t)=x+t$ for all $t\in \T_n$ and $x\in\field[R]^n$.
	Let $\{z^i\}$ and $\{u^i\}$ be coordinate systems on respectively
$\T_n$ and $\field[R]^n$ such that $z^i(t)=t^i$ and $u^i(x)=x^i$ for
 $t=(t^1,\dots,t^n)^\top\in \T_n$ and
 $x=(x^1,\dots,x^n)^\top\in \field[R]^n$.

    \begin{Prop}    \label{Prop4.2}
The left and right fundamental vector fields for $\T_n$ coincide and the
fundamental vector field associated with
$X=X^i\frac{\pd }{\pd z^i}\big|_{\bs{0}}\in T_{\bs{0}}(\T_n)$, $\bs{0}$ being
the zero vector of $\T_n$ (\ie of $\field[R]^n$), is
    \begin{equation}    \label{4.9}
\xi_X = X^i\frac{\pd }{\pd u^i} \in\Fields{X}(\field[R]^n) .
    \end{equation}
    \end{Prop}

    \begin{Proof}
	One can easily prove that $u^k\circ\lambda^x=u^k\circ\rho_x=x^k+z^k$
and
    \begin{equation}    \label{4.7}
(T_t(\lambda^x))\Bigl( \frac{\pd }{\pd z^i}\Big|_t \Bigr)
=
(T_t(\rho_x))\Bigl( \frac{\pd }{\pd z^i}\Big|_t \Bigr)
=
\frac{\pd }{\pd u^i}\Big|_{x+t} .
    \end{equation}
Therefore
    \begin{equation}    \label{4.8}
(T_t(\lambda^x))(Y_t)
=
(T_t(\lambda^x))(Y_t)
=
(Y_t)^i \frac{\pd }{\pd u^i}\Big|_{x+t}
    \end{equation}
for $Y_t=(Y_t)^i\frac{\pd }{\pd z^i}\big|_t\in T_t(\T_n)$. Putting here
$t=\bs{0}$ and $Y_{\bs{0}}=X$, we get~\eref{4.9}.
    \end{Proof}

	Obviously, the fundamental vector fields
    \begin{equation}    \label{4.10}
E_i = \frac{\pd }{\pd u^i} \in\Fields{X}(\field[R]^n)
    \end{equation}
form a basis for the set of fundamental fields of $\T_n$ in a sense that any
such field is their linear combination with \emph{constant} coefficients.

	The following result is almost trivial but nevertheless worth recording.

    \begin{Cor}    \label{Cor4.2}
A vector field on $\field[R]^n$ is with constant components (relative to
Cartesian coordinates) iff it is a fundamental vector field of $\T_n$.
    \end{Cor}



\subsection{Fundamental vector fields on \protect{$\field[R]^n$} induced
	    by \protect{$\GA(n,\field[R])$}}
	\label{Subsect4.3}

	The general affine group $\GA(n,\field[R])$ is a semidirect product
(sum) of the general linear group $\GL(n,\field[R])$ and the translation group
$\T_n$, $\GA(n,\field[R])=\GL(n,\field[R]) \rtimes \T_n$.~%
\footnote{~%
For a matrix realization of $\GA(n,\field[R])$ as a subgroup of
$\GL(n+1,\field[R])$, see~\cite[ch.~III, \S~3]{K&N-1}.%
}
	If $a_1,a_2\in\GL(n,\field[R])$ and $t_1,t_2\in \T_n$, the product of
the elements $a_1\rtimes t_1$ and $a_2\rtimes t_2$ in $\GA(n,\field[R])$ is
$(a_1\rtimes t_1)(a_2\rtimes t_2) := (a_1\cdot a_2)\rtimes(a_1\cdot t_2 + t_1)$.
	A left action
$\lambda \colon \GA(n,\field[R])\times\field[R]^n\to \field[R]^n$ is given by
    \begin{equation}    \label{4.10-1}
\lambda(a\rtimes t,x) := a\cdot x + t
    \end{equation}
for all $a=[a_i^j]\in\GL(n,\field[R])$, $t=(t^1,\dots,t^n)^\top\in \T_n$ and
$x\in\field[R]^n$.

    \begin{Prop}    \label{Prop4.3}
The fundamentals vector field $\xi_X$ on $\field[R]^n$ corresponding to
\(
X=
  X_j^i\frac{\pd }{\pd y_j^i}\big|_{\openone\rtimes \bs{0}}
+   X^i\frac{\pd }{\pd z^i}\big|_{\openone\rtimes \bs{0}}
\in T_{\openone\rtimes \bs{0}}(\GA(n,\field[R]))
\)
and the left action $\lambda$, given by~\eref{4.10-1}, is
    \begin{equation}    \label{4.12}
\xi_X
=
(X_j^i u^j + X^i) \frac{\pd }{\pd u^i} \in\Fields{X}(\field[R]^n)
    \end{equation}
in the coordinate systems $\{y_j^i,z^k\}$ and $\{u^i\}$ defined on the next
lines.
    \end{Prop}

    \begin{Proof}
	Define global coordinate systems $\{y_j^i,z^k\}$ on $\GA(n,\field[R])$
and $\{u^i\}$ on $\field[R]^n$ by
 $y_j^i(a\rtimes t):=a_j^i$, $z^i(a\rtimes t):=t^i$ and $u^k(x)=x^k$.
Then
$u^k\circ\lambda^x = u^l(x) y_l^k + z^k$ and consequently
    \begin{equation}    \label{4.11}
(T_{a\rtimes t}(\lambda^x))
\Bigl(
  Y_j^i\frac{\pd }{\pd y_j^i}\Big|_{a\rtimes t}
+   Z^i\frac{\pd }{\pd z^i}\Big|_{a\rtimes t}
\Bigr)
=
(Y_j^i u^j(x) + Z^i) \frac{\pd }{\pd u^i}\Big|_{a\cdot x+t}
\in T_{a\cdot x+t}(\field[R]^n)
    \end{equation}
for all $Y_j^i,Z^i\in\field[R]$. The assertion now follows from here and
definition~\ref{Defn4.1}.
    \end{Proof}

    \begin{Cor}        \label{Cor4.3}
A vector field on $\field[R]^n$ is an affine vector field (relative to some
Cartesian coordinates) iff it is a fundamental vector field of
$\GA(n,\field[R])$ represented on $\field[R]^n$ via the left action $\lambda$
described above.
    \end{Cor}



\subsection{Local left actions of \protect{$\GA(n,\field)$} on a manifold}
	\label{Subsect4.4}

The general affine group $\GA(n,\field)$  and its subgroups $\GL(n,\field)$
and $\T_n$ have natural (local) left actions on an arbitrary manifold $M$ of
dimension $\dim M=n$.

	Let $\lambda$ be the left action described in
subsection~\ref{Subsect4.3} with $\field[R]$ replaced with $\field$, $x\in M$
and $(V,v)$ be a chart of $M$ with $x$ in its domain, $x\in V$, and coordinate
diffeomorphism $v \colon V\to \field^n$. Define a mapping
 $L \colon \GA(n,\field)\times V\to V$ by
    \begin{equation}    \label{4.14}
L(a\rtimes t,x)
:= v^{-1}(a\cdot v(x)+t)
 = v^{-1}\circ\lambda(a\rtimes t,v(x))
    \end{equation}
for all $a\rtimes t\in\GA(n,\field)$, where $v(x)\in\field^n$ is considered as
a vector\ndash column. We have $L(a\rtimes t,x)\in V$ as $v(V)=\field^n$. Since
    \begin{equation}    \label{4.15}
L^x := L(\cdot,x) = v^{-1}\circ\lambda^{v(x)} \quad
L_{a\rtimes t} := L(a\rtimes t,\cdot) = v^{-1}\circ\lambda_{a\rtimes t}\circ v
    \end{equation}
and $\lambda \colon \GA(n,\field)\times\field^n\to \field^n$ is a left action
(see~\eref{4.10-1}), the mapping $L$ is a \emph{local} left action of
$\GA(n,\field)$ on $V\subseteq M$, which, obviously, depends on the chart
$(V,v)$.

	The (local) fundamental vector field  $\xi_X\in\Fields{X}(M)$
corresponding to $X\in T_{\openone\rtimes\bs{0}}(\GA(n,\field))$ and $L$ is
describe by the following proposition.
    \begin{Prop}    \label{Prop4.4}
The local fundamental vector field  $\xi_X\in\Fields{X}(M)$ corresponding to
\(
X=
  X_j^i\frac{\pd }{\pd y_j^i}\big|_{\openone\rtimes \bs{0}}
+   X^i\frac{\pd }{\pd z^i}\big|_{\openone\rtimes \bs{0}}
\in T_{\openone\rtimes \bs{0}}(\GA(n,\field[R])) ,
\)
where $\{y_i^j,z^k\}$ are the coordinates on $\GA(n,\field)$ introduced in
subsection~\ref{Subsect4.3},
and $L$ is
    \begin{equation}    \label{4.18}
\xi_X
=
(X_j^i v^j + X^i) \frac{\pd }{\pd v^i}
\in\Fields{X}(V)\subseteq\Fields{X}(M) .
    \end{equation}
    \end{Prop}

    \begin{Proof}
Let $\{u^i\}$ be the standard Cartesian coordinate system on $\field^n$ and
$\{v^i\}$ be the coordinate system defined by the chart $(V,v)$. Since
$u^i\circ v:=v^i$, we have
    \begin{equation}    \label{4.16}
v^k\circ L^x
= v^k\circ v^{-1}\circ\lambda^{v(x)}
= u^k\circ\lambda^{v(x)}
= u^l(v(x))y_l^k + z^k
= v^l(x)y_l^k + z^k .
    \end{equation}
Therefore
    \begin{equation}    \label{4.17}
(T_{a\rtimes t}(L^x))
\Bigl(
  Y_j^i\frac{\pd }{\pd y_j^i}\Big|_{a\rtimes t}
+   Z^i\frac{\pd }{\pd z^i}\Big|_{a\rtimes t}
\Bigr)
=
(Y_j^i v^j(x) + Z^i) \frac{\pd }{\pd v^i}\Big|_{L(a\rtimes t,x)}
\in T_{L(a\rtimes t,x)}(M)
    \end{equation}
for all $Y_j^i,Z^i\in\field$. The assertion now follows from here and
definition~\ref{Defn4.1}.
    \end{Proof}

	The following corollary is evident.

    \begin{Cor}    \label{Cor4.4}
A vector field on a $C^1$ (real or complex) manifold $M$ is an affine vector
field relative to a chart $(V,v)$ if and only if it reduces on $V$ to a
fundamental vector field of the general affine group $\GA(n,\field)$,
$n=\dim M$, represented on $M$, precisely on $V$, via the left action $L$
defined by~\eref{4.14}.
    \end{Cor}

	Since the groups $\GL(n,\field)$ and $\T_n$ are subgroups of
$\GA(n,\field)$, the above considerations can be applied  \emph{mutatis
mutandis} to them. Without going into details, this can be done as follows.

	The local left action of $\GL(n,\field)$ and $\T_n$ on $M$ in $(V,v)$
are given respectively by (cf.~\eref{4.14})
    \begin{subequations}    \label{4.20}
    \begin{align}    \label{4.20a}
L \colon \GL(n,\field)\times V & \to V  \colon (a,x)\mapsto  v^{-1}(a\cdot v(x))
\\		    \label{4.20b}
L \colon \T_n\times V & \to V  \colon (t,x)\mapsto v^{-1}(v(x) + t) .
    \end{align}
    \end{subequations}

    \begin{Prop}    \label{Prop4.5}
The fundamental vector fields corresponding to
 $X=X_j^i\frac{\pd }{\pd y_j^i} \in T_{\openone}(\GL(n,\field))$ and
 $X=X^i\frac{\pd }{\pd u^i} \in T_{\bs{0}}(\T_n)$ and the above actions $L$ are
respectively
    \begin{subequations}    \label{4.21}
    \begin{align}    \label{4.21a}
\xi_X & = X_j^i v^j\frac{\pd }{\pd v^i} \in \Fields{X}(V)\subseteq\Fields{X}(M)
\\		    \label{4.21b}
\xi_X & = X^i	 \frac{\pd }{\pd v^j} \in \Fields{X}(V)\subseteq\Fields{X}(M) .
    \end{align}
    \end{subequations}
    \end{Prop}

    \begin{Proof}
Since the mapping $a\mapsto a\rtimes\bs{0}$
(resp.\ $t\mapsto \openone\rtimes t$) realizes a homeomorphism from
$\GL(n,\field)$ (resp.\ $\T_n$) on $\GA(n,\field)$, we can assert that instead
of~\eref{4.17} now we have the equations
    \begin{subequations}    \label{4.22}
    \begin{align}    \label{4.22a}
(T_{a}(L^x))
\Bigl(   Y_j^i\frac{\pd }{\pd y_j^i}\Big|_{a} \Bigr)
& =
Y_j^i v^j(x) \frac{\pd }{\pd v^i}\Big|_{L(a,x)}
\in T_{L(a,x)}(M)
\\		    \label{4.22b}
(T_{t}(L^x))
\Bigl(   Z^i\frac{\pd }{\pd u^i}\Big|_{a} \Bigr)
& =
Z^i v^j(x) \frac{\pd }{\pd v^i}\Big|_{L(t,x)}
\in T_{L(t,x)}(M)
    \end{align}
    \end{subequations}
for respectively the left actions~\eref{4.20a} and~\eref{4.20b}. The
equations~\eref{4.21} follow immediately from here and
definition~\ref{Defn4.1}.
    \end{Proof}

	Now we shall record the evident analogue of corollary~\ref{Cor4.4}.

    \begin{Cor}    \label{Cor4.5}
A vector field on $C^1$ manifold $M$ is a linear (resp.\ constant) vector field
relative to a chart $(V,v)$ of $M$ if and only if it reduces on $V$ to a
fundamental vector field of the general linear (resp.\ translation) group
$\GL(n,\field)$ (resp.\ $\T_n$), $n=\dim M$, represented on $V$ via the left
action~\eref{4.20a} (resp.~\eref{4.20b}).
    \end{Cor}

	At last, we notice that the results of
subsections~\ref{Subsect4.1}--\ref{Subsect4.3} are special cases of the above
ones when $M=\field[R]^n$ and $(V,v)=(\field[R]^m,u)$, with $\{u^i\}$ being the
standard Cartesian coordinate system on $\field[R]^n$.


\subsection{Examples of ``non-standard'' actions of \protect{$\T_n$} and
	   \protect{$\GL(n,\field[R])$}}
	\label{Subsect4.5}

	The fundamental vector fields generally depend on the concrete action
of the group considered as it is clear from definition~\ref{Defn4.1}. The
purpose of the following lines is to illustrate this fact as well as an exception
of it.

	Consider the (left/right) action
    \begin{equation}    \label{4.23}
\lambda \colon \T_n\times\field[R]^n\to \field[R]^n
 \colon (t,x)\mapsto \lambda(t,x) = x\, \e^{s\cdot t}
    \end{equation}
of the translation group $\T_n$ on $\field[R]^n$. Here
$s=(s_1,\dots,s_n)^\top\in\field[R]^n$ is a fixed element and
$s\cdot t:=\sum_i s_it^i=s_it^i$ is the Euclidean product of $s$ and $t$. Using
the notation of subsection~\ref{Subsect4.2}, one finds that the fundamental
vector field corresponding to
$X= x^i\frac{\pd }{\pd z^i}\in T_{\bs{0}}(\T_n)$
and $\lambda$ is
    \begin{equation}    \label{4.24}
\xi_X = (X^i s_i) x^k\frac{\pd }{\pd u^k} \in \Fields{X}(\field[R]^n).
    \end{equation}
Therefore linear vector fields of the type $cx^k\frac{\pd }{\pd u^k}$, $c$ being
a real constant, are fundamental vector fields of $\T_n$ and $\lambda$ for
suitable choice of $X$ and/or $s$ and \emph{vice versa}. Therefore the
fundamental vector fields of $\T_n$ relative to the representation~\eref{4.23}
are linear vector fields while the ones relative to representation via
translations are constant vector fields.

	As a second example, let us investigate the (left) action
    \begin{equation}    \label{4.25}
\lambda \colon \GL(n,\field[R]) \times\field[R]^n\to \field[R]^n
 \colon (a,x)\mapsto \lambda(t,x) = ax (\det a)^q ,
    \end{equation}
for some number $q\in\field[N]\cup\{0\}$,
of $\GL(n,\field[R])$ on $\field[R]^n$; the case $q=0$ being the one
considered in subsection~\ref{Subsect4.1}. Calculating the appearing
in~\eref{4.0-2} derivatives
(notice that $\frac{\pd \det a}{\pd a_i^j}\big|_{a=\openone} = \delta_j^i $
with $\delta_j^i$ being the Kronecker deltas),
we see
that the fundamental vector field corresponding to $\GL(n,\field[R])$, the
action~\eref{4.25} and a vector $X\in T_{\openone}(\GL(n,\field[R]))$ is
    \begin{equation}    \label{4.26}
\xi_X
     = (X_j^i x^j + q X_j^j x^i) \frac{\pd }{\pd u^i}
     = (X_j^i + q X_k^k \delta_j^i) x^j \frac{\pd }{\pd u^i} .
    \end{equation}
The set of these fundamental vector fields coincides with the one of linear
vector fields as, if
$C_j^i=X_j^i + q X_k^k \delta_j^i$, then
$X_j^i=C_j^i - \frac{q}{1+qn} C_k^k\delta_j^i$. Hence it coincides with the
one of $\GL(n,\field[R])$ represented on $\field[R]^n$ via left
multiplication. However, the particular fundamental vector fields corresponding
to a concrete vector $X\in T_{\openone}(\GL(n,\field[R]))$ are different for
the two actions considered (unless $q=0$ when they are identical).



\section {Conclusion}
\label{Conclusion}

      As we already said in section~\ref{Introduction}, this paper reviews the
concepts of affine and fundamental vector fields on a manifold. The main
conclusions from it are that the affine, linear and constant vector fields on a
manifold are in a bijective correspondence with the fundamental vector fields
on it of respectively general affine, general linear and translation groups
(locally) represented on the manifold via the described in this work left
actions; in a case of the manifold $\field^n=\field[R]^n,\field[C]^n$, the
actions mentioned have the usual meaning of affine, linear and translation
transformations.

	Equations~\eref{4.0-2} can serve as a ground for studding a problem
inverse to the one of finding fundamental vector fields, \viz to be found a Lie
group and its action on a manifold if some set of vector fields plays a role
of set of its fundamental vector fields.  Precisely, given numbers
$X^\mu\in\field$, $\mu=1,\dots,N\in\field[N]$, and vector fields~%
\footnote{~%
Only vector fields of the type~\eref{Con.1} can be fundamental vector
fields of some group according to~\eref{4.0-2}.
}
    \begin{equation}    \label{Con.1}
\xi_X = X^\mu f_\mu^k \frac{\pd }{\pd u^k}
    \end{equation}
for some functions $f_\mu^k$ on $M$. Does there exists a Lie group $G$ with
$\dim G=N$ and a left/right action of $G$ on $M$ for which $\xi_X$ is the
fundamental vector field corresponding to
 $X=X^\mu\frac{\pd }{\pd y^m}\in T_eG$? In particular, the group $G$ can be
given and one should look for the existence/non-existence of (one or more)
actions with the last property.

	Comparing~\eref{Con.1} and, e.g.,~\eref{4.0-2a}, we get
    \begin{equation}    \label{Con.2}
u^k\circ\lambda(a,x)
=
u^k(x) + f_\mu^k(x) [y^\mu(a)-y^\mu(e)]
 + f_{\mu\nu}^k(a,x)[y^\mu(a)-y^\mu(e)][y^\nu(a)-y^\nu(e)]
    \end{equation}
due to $\lambda(e,x)\equiv x$. Here $f_{\mu\nu}^k$ and their first partial
derivatives relative to $y^\mu$ are bounded functions. To define a left action
$\lambda$ via this equation one needs to ensure that
$\lambda_{ab}=\lambda_a\circ\lambda_b$ for all $a,b\in G$. As a result
of~\eref{Con.2}, this requirement is equivalent to
    \begin{multline}    \label{Con.3}
f_\mu^k(x) [y^\mu(ab)-y^\mu(e)]
 + f_{\mu\nu}^k(ab,x)[y^\mu(ab)-y^\mu(e)][y^\nu(ab)-y^\nu(e)]
\\ =
f_\mu^k(x) [y^\mu(b)-y^\mu(e)]
 + f_{\mu\nu}^k(b,x)[y^\mu(b)-y^\mu(e)][y^\nu(b)-y^\nu(e)]
\\ +
f_\mu^k(\lambda(b,x)) [y^\mu(a)-y^\mu(e)]
+
f_{\mu\nu}^k(a,\lambda(b,x))[y^\mu(a)-y^\mu(e)][y^\nu(a)-y^\nu(e)]
    \end{multline}
which should hold for all $a,b\in G$ and $x\in M$. This is a system of
equations for the left action $\lambda$ (involved directly and via its
expansion~\eref{Con.2}) and the for the multiplication
$G\times G\ni(a,b)\mapsto ab\in G$ in $G$.


\section*{Acknowledgments}

	The author's interest in the problems, considered in the present paper,
arose from a terminological discussion with Prof.~ Maido Rahula (Institute of
pure mathematics, Faculty of mathematics and informatics, University of Tartu,
Tartu, Estonia) about should a fundamental vector field be called in this way or
to be termed `operator of a group' (or `group operator'). Respectively, this
work is partially done within the Joint research project ``Vector fields and
symmetries'' within the bilateral agreement between the Bulgarian academy of
sciences and the Estonian academy of sciences.

	This work was partially supported by the National Science Fund of
Bulgaria under Grant No.~F~1515/2005.


\addcontentsline{toc}{section}{References}
\bibliography{bozhopub,bozhoref}

\begin{thebibliography}{10}

\bibitem{Olver/LieGroups}
P.~J. Olver.
\newblock {\em Applications of Lie groups to differential equations}.
\newblock Springer-Verlag, Yew York, 2 edition, 1993.
\newblock 513~p. First ed.: Springer-Verlag, New York, 1986. Russian
  translation: Mir, Moscow, 1989.

\bibitem{Warner}
F.~W. Warner.
\newblock {\em Foundations of differentiable manifolds and Lie groups}.
\newblock Springer-Verlag, New York-Berlin-Heidelberg-Tokyo, 1983.
\newblock Russian translation: Mir, Moscow, 1987.

\bibitem{K&N-1}
S.~Kobayashi and K.~Nomizu.
\newblock {\em Foundations of Differential Geometry}, volume~I.
\newblock Interscience Publishers, New York-London, 1963.
\newblock Russian translation: Nauka, Moscow, 1981.

\bibitem{Matveev1974}
N.~M. Matveev.
\newblock {\em Methods of integrating ordinary differential equations}.
\newblock High school, Minsk, 1974.
\newblock In Russian.

\bibitem{KMS-1993}
Ivan Kol\'a\v{r}, Peter~W. Michor, and Jan Slov\'ak.
\newblock {\em Natural operations in differential geometry}.
\newblock Springer, Berlin - Heidelberg, 1993.

\bibitem{Poor}
Walter~A. Poor.
\newblock {\em Differential geometric structures}.
\newblock McGraw-Hill Book Company Inc., New York, 1981.

\bibitem{Greub&et_al.-2}
W.~Greub, S.~Halperin, and R.~Vanstone.
\newblock {\em Lie groups, principal bundles, and characteristic classes},
  volume~2 of {\em Connections, Curvature, and Cohomology}.
\newblock Academic Press, New York and London, 1973.

\bibitem{Chebotarev-1940}
N.~G. Chebotarev.
\newblock {\em Theory of the Lie groups}.
\newblock GITTL, Moscow, 1940.
\newblock In Russian.

\bibitem{Lie/Groups}
S.~Lie.
\newblock {\em Theorie der Transformationsgruppen}, volume I, II, III.
\newblock \, Leipzig, 1888, 1890, 1893.

\bibitem{Cantwell}
Brian~J. Cantwell.
\newblock {\em Introduction to symmetry analysis}.
\newblock Cambridge University Press, 2002.

\end{thebibliography}
\bibliographystyle{unsrt}
\addcontentsline{toc}{subsubsection}{This article ends at page}



\end{document}